\documentclass[]{amsart}

\usepackage{amssymb,latexsym,mathrsfs}
\numberwithin{equation}{section}

\newtheorem{thm}{Theorem}[section]
\newtheorem*{thm*}{Theorem}

\newtheorem{conjecture}[thm]{Conjecture}

\theoremstyle{remark}

\newcommand{\Ker}{\mathrm{Ker}}
\newcommand{\GL}{\mathrm{GL}}
\newcommand{\Li}{\mathrm{Li}}
\newcommand{\Lscr}{\mathscr{L}}
\newcommand{\Cscr}{\mathscr{C}}
\newcommand{\Fscr}{\mathscr{F}}
\newcommand{\Iscr}{\mathscr{I}}
\newcommand{\Cl}{\mathrm{Cl}}
\newcommand{\Cbf}{\mathbf{C}}
\newcommand{\Dbf}{\mathbf{D}}

\newcommand{\ep}{\varepsilon}

\newcommand{\con}{\equiv}

\newcommand{\ndiv}{\nmid}
\newcommand{\modd}[1]{\; ( \text{mod} \; #1)}
\newcommand{\bstack}[2]{#1 \atop #2}

\newcommand{\maps}{\rightarrow}
\newcommand{\intersect}{\cap}

\newcommand{\Gal}{\mathrm{Gal}}

\newcommand{\al}{\alpha}
\newcommand{\be}{\beta}

\newcommand{\ga}{\gamma}
\newcommand{\del}{\delta}
\newcommand{\Del}{\Delta}

\newcommand{\sig}{\sigma}

\newcommand{\Ocal}{\mathcal{O}}

\newcommand{\N}{\mathbb{N}}
\newcommand{\Q}{\mathbb{Q}}

\newcommand{\Z}{\mathbb{Z}}

\newcommand{\afr}{\mathfrak{a}}
\newcommand{\bfr}{\mathfrak{b}}

\newcommand{\pfr}{\mathfrak{p}}

\newcommand{\Nfr}{\mathfrak{N}}

\newcommand{\beq}{\begin{equation}}
\newcommand{\eeq}{\end{equation}}

\begin{document}

\title[Counting problems: class groups, primes, and number fields]{Counting problems:\\ class groups, primes, and number fields}

\author{Lillian B. Pierce}
\address[Lillian B. Pierce]{Mathematics Department, Duke University, Durham NC 27708 USA} \email{pierce@math.duke.edu}

\begin{abstract}
Each number field has an associated finite abelian group, the class group, that records certain properties of arithmetic within the ring of integers of the field. 
The class group is well-studied, yet also still mysterious. A central conjecture of Brumer and Silverman states that for each prime $\ell$, every number field has the property that its class group has very few elements of order $\ell$, where ``very few'' is measured relative to the absolute discriminant of the field. This paper surveys recent progress toward this conjecture, and outlines its close connections to counting prime numbers, counting number fields of fixed discriminant, and counting number fields of bounded discriminant.
 
\end{abstract}

 %\dedicatory{dedicated to my teachers and collaborators}
\maketitle

%\setcounter{tocdepth}{2}
% \tableofcontents

\section{Historical prelude}
In a 1640 letter to Mersenne, Fermat 
stated that an odd prime $p$ satisfies $p=x^2 +y^2$ if and only if $p \con 1 \modd{4}$. 
Roughly 90 years later,
Euler learned of Fermat's statement through correspondence with Goldbach, and by 1749, he worked out a proof.
This   fits into a bigger question, which Euler studied as well: for each $n\geq 1$, which primes can be written as $p=x^2 +ny^2?$ Even more generally: which binary quadratic forms $ax^2 + bxy + cy^2$ represent a given integer $m$? This question also motivated work of Lagrange   and Legendre, and then appeared in Gauss's celebrated 1801 work  \emph{Disquisitiones Arithmeticae}; see \cite{Cox13}.

Gauss partitioned binary quadratic forms of discriminant $D = b^2 - 4ac$ into equivalence classes under $\mathrm{SL}_2(\Z)$ changes of variable.
 (Here we will speak only of \emph{fundamental} discriminants $D$; for notes on the original setting, see \cite{Sta07}.)
Gauss showed that for each $D$ there are finitely many such classes (the cardinality is the \emph{class number}, denoted $h(D)$), and verified that the classes obey a group law (composition).
Based on extensive computation, Gauss noticed that as $D \maps -\infty$, small class numbers stopped appearing, writing: 
``\emph{Nullum dubium esse videtur, quin series adscriptae
revera abruptae sint...Demonstrationes autem}
rigorosae \emph{harum observationum perdifficiles esse videntur.}''
(``It seems beyond doubt that the sequences written down do indeed break off...However, \emph{rigorous}
proofs of these observations appear to be most difficult'' \cite[p. 13]{GSS07}.) As $D \maps + \infty$, quite different behavior seemed to appear, leading to a conjecture that $h(D)=1$ for infinitely many $D>0$.  

It is hard to exaggerate the interest these two conjectures have generated.
In the 1830's, Dirichlet proved a class number formula, relating the class number $h(D)$  of a (fundamental) discriminant $D$ to the value $L(1,\chi)$ of an $L$-function associated to a real primitive character $\chi$ modulo $D$. Consequently, throughout the 1900's, Gauss's questions were studied via the theory of the complex-variable functions $L(s,\chi)$.
A remarkable series of works by Hecke, Deuring, Mordell, and Heilbronn confirmed that for $D<0$ the class number $h(D)$ attains any value only finitely many times. How many times?  Famously, work of  Heegner, Baker, and Stark proved there are 9 (fundamental) discriminants $D<0$ with class number $1.$ In full generality, Goldfeld showed an effective lower bound for $h(D)$ when $D<0$ would follow from a specific case of the Birch--Swinnerton-Dyer conjecture, which was then verified by Gross and Zagier; see \cite{Gol85}. Now, for each $1 \leq N \leq 100$, one may find the number of discriminants $D<0$ with $h(D)=N$   in \cite{Wat04}. 
As for the other conjecture, that infinitely many (fundamental) discriminants $D>0$ have class number 1, this remains open, and very mysterious.
These historical antecedents   hint at the intertwined currents of ``counting'' and the analytic study of $L$-functions, which will also be present in the work we will survey.

We briefly mention another historical motivation for the study of class numbers: the failure of unique factorization.
For example, in the ring $\Z[\sqrt{-5}]$,  $21 = 3 \cdot 7$ but also factors into irreducible, non-associated factors as $(1+2\sqrt{-5})(1 - 2\sqrt{-5})
$.  
Here is a problem where the failure of unique factorization has an impact. Suppose one is searching for nontrivial solutions $x,y,z\in \N$ to the equation 
$ x^p + y^p = z^p $
for a prime $p \geq 3.$
If  a nontrivial solution $(x,y,z)$ exists, then for $\zeta_p$  a $p$-th root of unity, we could write
$
y \cdot y \cdots y = (z- x)(z-\zeta_p x) \cdots (z - \zeta_p^{p-1}x).
$
If $\Z[\zeta_p]$ possesses   unique factorization, two such  factorizations cannot exist, so $(x,y,z)$ cannot exist---verifying  Fermat's Last Theorem  for this exponent $p$. 
But to the disappointment of many, unique factorization fails in $\Z[\zeta_p]$  for infinitely many $p$. As Neukirch writes, ``Realizing the failure of unique factorization in general has led to one of the grand events in the history of number theory, the discovery of ideal theory by Eduard Kummer''   \cite[Ch. I \S 3]{Neu}.

\subsection{The class group}
Let $K/\Q$ be a number field of degree $n$,  with associated ring of integers $\Ocal_K$.
Every proper integral ideal $\afr \subset \Ocal_K$ factors into a product of prime ideals $\pfr_1 \cdots \pfr_k$ in a unique way (salvaging the notion of unique factorization).  
Moreover, the fractional ideals of $K$ form an abelian group $J_K$, the free abelian group on the set of nonzero prime ideals of $\Ocal_K$. 
In the case that every ideal in $J_K$ belongs to the subgroup $P_K$ of principal ideals, then $\Ocal_K$ is a principal ideal domain, and unique factorization holds in $\Ocal_K$. But more typically, some ``expansion'' occurs when passing to ideals; the class group of $K$ is defined to measure this.

  The class group of $K$ is
 the quotient group
\[\Cl_K = J_K/P_K.\]
The elements in $\Cl_K$ are ideal classes, and
the cardinality $|\Cl_K|$ is  the class number. The quotient $J_K/P_K$ is trivial (so that every ideal is a principal ideal, and $|\Cl_K|=1$) precisely when unique factorization holds in  $\Ocal_K$. 
(Thus  the above strategy    for Fermat's Last Theorem  works for $p$ if $|\Cl_{\Q(\zeta_p)}|=1$.
In fact, Kummer showed that as long as the class number of $\Q(\zeta_p)$ is indivisible by $p$, the argument can be salvaged; see \cite{Edw96}. Such a prime is called a regular prime. Here is an open question: are there infinitely many regular prime numbers?)

By a result of Minkowski in the geometry of numbers, every ideal class in $\Cl_K$ contains an integral ideal $\bfr$ with norm $\Nfr(\bfr) = (\Ocal_K : \bfr)$ satisfying 
\beq\label{Norm_b}
\Nfr(\bfr) \leq (2/\pi)^s \sqrt{D_K},
\eeq
where $D_K = |\mathrm{Disc}(K/\Q)|$ and $s$ counts the pairs of complex embeddings of $K$. As there are finitely many integral ideals of any given norm, Landau deduced (see \cite[Thm. 4.4]{Nar90}):
\beq\label{Minkowski_bound}
|\Cl_K| \ll_n D_K^{1/2} \log^{n-1}D_K.
\eeq
 In particular, the class group of a number field $K$ is always a finite abelian group.
  (Throughout, $A \ll_\kappa B$ indicates that there exists a constant $C_\kappa$ such that $|A| \leq C_\kappa B$.)

When $K=\Q(\sqrt{D})$ is a  quadratic field, this relates in a precise way to Gauss's construction of the class number for binary quadratic forms of discriminant $D$ (see \cite{Bha06ICM}).
In modern terms, Gauss asked whether for each $h \in \N$, there are finitely many imaginary quadratic fields $K$ with $|\Cl_K|=h?$ (Yes.) 
Are there infinitely many real quadratic fields $K$ with $|\Cl_K|=1$? (We don't know.) In fact, here  is an open question: are there infinitely many number fields, of arbitrary degrees, with class number 1? 
Here is another open question: are there infinitely many number fields, of arbitrary degrees, with  bounded class number?
These difficult questions must consider the regulator $R_K$ of the field $K$, due to the (ineffective) inequalities by Siegel (for quadratic fields) and Brauer (in general) \cite[Ch. 8]{Nar90}:
\[
D_K^{1/2-\ep} \ll_{n,\ep} |\Cl_K| R_K \ll_{n,\ep} D_K^{1/2+\ep}, \quad \text{for all $\ep>0$.}
\]

\section{The $\ell$-torsion Conjecture}
 In addition to studying the size of the class group, it is also natural to study its  structure. We will focus on the  $\ell$-torsion subgroup, defined for each integer $\ell \geq 2$  by
\[ \Cl_K[\ell] = \{ [\afr] \in \Cl_K :   [\afr]^\ell = \mathrm{Id}\}.
\]
For example, the class number is divisible by a prime $\ell$ precisely when $|\Cl_K[\ell]|>1$.
Related problems include: to study the exponent of the class group, or to count how many number fields of a certain degree have class number divisible, or indivisible, by a given prime $\ell$. Such problems are addressed for imaginary quadratic fields in \cite{Sou00,HB07a,HB08,Bec17}.

In this survey, we will focus on upper bounds for the $\ell$-torsion subgroup.
The Minkowski bound (\ref{Minkowski_bound}) provides an upper bound for any field of degree $n$, and all $\ell$: 
\beq\label{class_trivial}
1 \leq |\Cl_K[\ell]| \leq |\Cl_K| \ll_{n,\ep} D_K^{1/2+\ep}, \qquad \text{for all $\ep>0$.}
\eeq
Our subject is a conjecture on the size of the $\ell$-torsion subgroup, which suggests that (\ref{class_trivial}) is far from the truth.
  We will focus primarily on cases when $\ell$ is prime, since $|\Cl_K[m]|$ is multiplicative as a function of $m$, and for a prime $\ell$, $|\Cl_K[\ell^t]| \leq |\Cl_K[\ell]|^t$ (see \cite{PTBW21}).

\begin{conjecture}[$\ell$-torsion Conjecture]\label{conj_class}
 Fix a degree $n \geq 2$ and a prime $\ell$. Every number field $K/\Q$ of degree $n$ satisfies  
$
|\Cl_{K}[\ell]| \ll_{n,\ell,\ep}D_K^\ep$ for all $\ep>0$.
\end{conjecture}

This conjecture is due to Brumer and Silverman,  in the more precise form: is it always true that
$ \log_\ell |\Cl_K[\ell]|\ll_{n,\ell} \log D_K/ \log \log D_K$ \cite[Question {$\Cl(\ell,d)$}]{BruSil96}? 
Brumer and Silverman were motivated by counting elliptic curves of fixed conductor. Subsequently, this conjecture has appeared in many further contexts, including:
  bounding   the ranks of elliptic curves \cite[\S 1.2]{EllVen07}; bounding Selmer groups and  ranks of   hyperelliptic curves \cite{BSTTTZ20};  counting number fields   \cite[p. 166]{Duk98}; equidistribution of CM points on Shimura varieties  \cite[Conjecture 3.5]{Zha05};   counting non-uniform lattices in semisimple Lie groups \cite[Thm. 7.5]{BelLub19}.

Conjecture \ref{conj_class} is known to be true for the degree $n=2$ and the prime $\ell=2$, when it follows from the genus theory of Gauss (see \cite[Ch. 8.3]{Nar90}). This is the only case in which it is known.
Nevertheless, starting in the early 2000's, significant progress has been made. The purpose of this survey is to give some insight into the wide variety of methods developed in recent   work toward the conjecture. 
As an initial measure of progress, we define:\\
\textbf{Property $\Cbf_{n,\ell}(\Del)$:}
\emph{Fix a degree $n \geq 2$ and a prime $\ell$. Property $\Cbf_{n,\ell}( \Del)$ holds if  for all number fields $K/\Q$ of degree $n$, 
 $ |\Cl_K[\ell]| \ll_{n,\ell,\Del,\ep} D_K^{\Del + \ep}$ for all $\ep>0$.
  }\\
 Gauss proved that $\Cbf_{2,2}(0)$ holds.
  Until recently, no other case with $\Del<1/2$ was known.

The first   progress  was for imaginary quadratic fields. Suppose $K = \Q(\sqrt{-d})$ for a square-free integer $d>1$, and suppose that $[\afr]$ is a nontrivial element in $\Cl_K[\ell]$ for a prime $\ell \geq 3;$ thus $[\afr]^\ell$ is the principal ideal class. Then by the Minkowski bound (\ref{Norm_b}), there exists an integral ideal $\bfr$ in $[\afr]$ such that $\Nfr(\bfr) \ll d^{1/2}.$ Moreover $\bfr^\ell$ is principal, say generated by $(y+z \sqrt{-d})/2$ for some integers $y,z$,
and so
$(\Nfr(\bfr))^\ell=\Nfr(\bfr^\ell) =  (y^2 + dz^2)/4$. Consequently,   $|\Cl_K[\ell]|$ can be dominated (up to a factor $d^\ep$) by the number of integral solutions to
 \beq\label{surface_region}
 4x^\ell = y^2 + dz^2, \qquad \text{with $x\ll d^{1/2},\; y \ll d^{\ell/4}, \; z \ll d^{\ell/4 - 1/2}$}.
 \eeq
 When $\ell=3$, this can be interpreted in several ways: counting solutions to a congruence $y^2 = 4x^3 \modd{d}$; counting perfect square values of the polynomial $f(x,z)= 4x^3 - dz^2$; or counting integral points on a family of Mordell elliptic curves $y^2 = 4x^3 - D,$ with $D=dz^2.$ 
Pierce used the first two perspectives, and Helfgott and Venkatesh used the third perspective, to prove for the first time that property $\Cbf_{2,3}(\Del)$ holds for some $\Del<1/2$  \cite{Pie05,Pie06,HelVen06}. (The Scholz reflection principle shows that $\log_3 |\Cl_{\Q(\sqrt{-d})}[3]|$ and $\log_3 |\Cl_{\Q(\sqrt{3d})}[3]|$ differ by at most 1, so results for $3$-torsion apply comparably to both real and imaginary quadratic fields \cite{Sch32}.)
  When $\ell \geq 5$,  the region in which $x,y,z$  lie in  (\ref{surface_region}) becomes inconveniently large relative to the trivial bound (\ref{class_trivial}).
 Here is an open question: for a prime $\ell \geq 5$, are there at most $\ll d^{\Del}$ integral solutions to (\ref{surface_region}), for some $\Del<1/2$?

Recently, Bhargava, Taniguchi, Thorne, Tsimerman, and Zhao  made a breakthrough on property $\Cbf_{n,2}(\Del)$ for all $n \geq 3$. Fix a prime $\ell$ and a number field $K$ of degree $n$. Given any nontrivial ideal class $[\afr] \in \Cl_K[\ell]$, they show it contains an integral ideal $\bfr$ with $\bfr^\ell$ a principal ideal generated by an element   $\be$ lying in a well-proportioned ``box.'' 
By an ingenious geometry of numbers argument, they show the number of such generators $\be$ in the box is $\ll D_K^{\ell/2-1/2}.$ If $\ell \geq 3$ this far exceeds the trivial bound (\ref{class_trivial}), but if $\ell=2$, it slightly improves it. The striking refinement comes  by recalling that any  $\beta$    of interest must also have $|N_{K/\Q}(\be)| = \Nfr(\bfr^\ell) = (\Nfr(\bfr))^\ell$ be a perfect $\ell$-th power of an integer, say $y^\ell$. For $\ell=2$ they apply a celebrated result of Bombieri and Pila to count  integral solutions $(x,y)$ to the degree $n$  equation $N_{K/\Q}(\be+x) = y^2$ \cite{BomPil89}.
This strategy proves property $\Cbf_{n,2}(1/2-1/2n)$ holds for all degrees $n \geq 3.$ Further refinements for degrees $3,4$ show $\Cbf_{3,2}(0.2785...)$ and $\Cbf_{4,2}(0.2785...)$ hold; see \cite{BSTTTZ20}.

Only two further nontrivial cases of property $\Cbf_{n,\ell}(\Del)$ are known, and for these we introduce the Ellenberg--Venkatesh criterion.

\subsection{The Ellenberg--Venkatesh criterion}
An important criterion for bounding $\ell$-torsion in the class group of a number field $K$  relies on counting small   primes that are non-inert in $K$.   The germ of the idea, which has been credited independently to Soundararajan and Michel, goes as follows. Suppose for example that $K = \Q(\sqrt{-d})$ is an imaginary quadratic field with $d$ square-free, and $\ell$ is an odd prime. Let $H$ denote  $\Cl_K[\ell]$. Then 
$|H| = |\Cl_K| / [\Cl_K:H],$
and to show that $|H|$ is small, it suffices to show that the index $[\Cl_K:H]$ is large. Now suppose that $p_1 \neq p_2$ are rational primes not dividing $2d$ that both split in $K$, say $p_1 = \pfr_1 \pfr_1^\sig$ and $p_2 = \pfr_2 \pfr_2^\sig,$ where $\sig$ is the nontrivial automorphism of $K.$ We claim that as long as $p_1,p_2$ are sufficiently small, then $\pfr_1$ and $\pfr_2$ must represent different cosets of $H$. For supposing to the contrary that $\pfr_1 H = \pfr_2 H$, one deduces that $\pfr_1 \pfr_2^\sig \in H$ so that $(\pfr_1 \pfr_2^\sig)^\ell$ is a principal ideal, say generated by $(y+z\sqrt{-d})/2$, for some $y,z \in \Z$. Taking norms shows
\beq\label{prime_equation}
4 (p_1p_2)^\ell = y^2 + dz^2.
\eeq
If   $p_1,p_2 < (1/4) d^{1/(2\ell)},$ this forces $z=0$, which yields a contradiction, since $4(p_1p_2)^\ell$ cannot be a perfect square. This proves the claim. In particular, if there are $M$ such distinct primes  $p_1,\ldots, p_M < (1/4) d^{1/2\ell}$ with $p_j \ndiv 2d$ and  $p_j$ split in $K$, then $|\Cl_K[\ell]| \leq |\Cl_K|M^{-1}.$

Ellenberg and Venkatesh significantly generalized this strategy to prove an influential criterion, which we state in the case of extensions of $\Q$ \cite{EllVen07}.  (Throughout this survey, we will  focus for simplicity on extensions of $\Q$, but many of the theorems and questions we mention have analogues in the literature over any fixed number field.) 
\newline
\textbf{Ellenberg--Venkatesh criterion:} \emph{Suppose $K/\Q$ is a number field  of degree $n \geq 2$, fix an integer $\ell\geq 2$, and fix $\eta< \frac{1}{2\ell(n-1)}.$ Suppose that there are $M$ prime ideals  $\pfr_1,\ldots,\pfr_M \subset \Ocal_K$ such that each $\pfr_j$ has norm $\Nfr (\pfr_j)<D_K^\eta$, $\pfr_j$ is unramified in $K$ and $\pfr_j$ is not an extension of a prime ideal from any proper subfield of $K$. Then
\beq\label{EV_criterion}
|\Cl_K[\ell]| \ll_{n,\ell,\ep} D_K^{\frac{1}{2}+ \ep}M^{-1},\qquad \text{for all $\ep>0$.}
\eeq
}
(A prime ideal $\pfr \subset \Ocal_K$ lying above a prime $p \in \Q$ is unramified in $K/\Q$ if $\pfr^2 \ndiv p\Ocal_K$; a prime  ideal $\pfr \subset \Ocal_K$ is an extension of a prime ideal in a proper subfield $K_0 \subset K$ if there exists a prime ideal $\pfr_0 \subset \Ocal_{K_0}$ such that $\pfr = \pfr_0 \Ocal_K.$)
For example,
if $p< D_K^\eta$ is a rational prime that splits completely in $K$, so that $p\Ocal_K = \pfr_1 \cdots \pfr_n$ for distinct prime ideals $\pfr_j$, then each $\pfr_j$ satisfies the hypotheses of the criterion.
In particular, if  $M$ rational (unramified) primes $p_1,\ldots, p_M < D_K^\eta$  split completely in $K$, then (\ref{EV_criterion}) holds. Alternatively, it suffices to exhibit prime ideals $\pfr_j \subset \Ocal_K$ of degree 1, since such a prime ideal cannot be an extension of a prime ideal from  a proper subfield.

Here is one of   Ellenberg and Venkatesh's striking   applications, which shows that   $\Cbf_{2,3}(1/3)$ holds---the current record for $n=2, \ell=3$.
Fix a large square-free integer $d>1$. Any prime $p \ndiv 6d$ that is inert in $\Q(\sqrt{-3})$ must split either in $\Q(\sqrt{d})$ or in $\Q(\sqrt{-3d}).$ 
Thus for any $\eta<1/6$, at least one field $K \in \{\Q(\sqrt{d}),\Q(\sqrt{-3d})\}$ has  a positive proportion of the primes $(1/2)d^{\eta} \leq p \leq d^{\eta}$ split in $K$. By the Ellenberg--Venkatesh criterion (\ref{EV_criterion}), this field $K$ then has the property that $|\Cl_K[3]| \ll D_K^{1/3+\ep}$ for all $\ep>0$. By the Scholz reflection principle, this bound also applies to the other field in the pair, and $\Cbf_{2,3}(1/3)$ holds.

The Scholz reflection principle has also been generalized by Ellenberg and Venkatesh to bound $\ell$-torsion (for odd primes $\ell$) in class groups  of even-degree extensions of certain number fields. 
In particular, by pairing their criterion with a reflection principle, they show that $\Cbf_{3,3}(1/3)$ holds and $\Cbf_{4,3}(\Del)$ holds for some $\Del<1/2$ \cite[Cor. 3.7]{EllVen07}. 
This concludes the list of degrees $n$ and primes $\ell$ for which property $\Cbf_{n,\ell}(\Del)$ is known   for some $\Del<1/2.$ 

Here are open problems: reduce the value $\Del<1/2$ for which $\Cbf_{n,\ell}(\Del)$ holds, when $n \geq 3$ and $\ell=2$, or when $n=2,3$ or $4$ and $\ell=3.$
For $n =2, 3$ or $4$ and a prime $\ell \geq 5$,  prove for the first time that $\Cbf_{n,\ell}(\Del)$ holds for some $\Del<1/2$.
For $n\geq 5$ and a prime $\ell \geq 3$,   prove for the first time that $\Cbf_{n,\ell}(\Del)$ holds for some $\Del<1/2.$

The Ellenberg--Venkatesh criterion underlies most of the significant recent progress on bounding $\ell$-torsion in class groups. What is the best result it can imply? Assuming the Generalized Riemann Hypothesis,  given any number field $K/\Q$ of degree $n$ with $D_K$ sufficiently large, a positive proportion of primes $p < D_K^{\eta}$ split completely in $K$, implying  
\beq\label{GRH_bound}
|\Cl_K[\ell]| \ll_{n,\ell,\ep} D_K^{\frac{1}{2} - \frac{1}{2\ell(n-1)} + \ep}, \qquad \text{for all $\ep>0.$}
\eeq
As this is a useful benchmark, we will call this the GRH-bound, and for convenience  set $\Del_{\mathrm{GRH}}=\frac{1}{2} - \frac{1}{2\ell(n-1)}$ once $n,\ell$ have been fixed. 
Thus if GRH is true, for each $n, \ell$, property $\Cbf_{n,\ell}(\Del_{\mathrm{GRH}})$ holds. There has been intense interest   in proving this without assuming GRH, and this will be our next topic.

 \section{Families of fields}
So far we have considered, for each degree $n$, the ``family'' of number fields $K/\Q$ of degree $n$. Let us formalize this, letting $\Fscr_n(X)$ be the set of all degree $n$ extensions $K$ of $\Q$, with $D_K = |\mathrm{Disc}(K/\Q)|\leq X$; let $\Fscr_n = \Fscr_n(\infty)$.
It is helpful at this point to consider more specific families of fields of a fixed degree.  
For example,  we could define $\Fscr_2^-(X)$    to be the  set of imaginary quadratic fields $K$ with $D_K \leq X$, and similarly $\Fscr_2^+(X)$ for real quadratic fields.
In general, given a transitive subgroup $G \subset S_n$, define the family
\beq\label{family_dfn}
  \Fscr_n(G;X) = \{ K/\Q:   \deg K/\Q = n, \Gal(\tilde{K}/\Q) \simeq G, D_K \leq X\},
\eeq
where all $K$ are in a fixed algebraic closure $\overline{\Q}$, $\tilde{K}$ is the
Galois closure of $K/\Q$, the Galois group is considered as a permutation
group on the $n$ embeddings of $K$ in $\overline{\Q}$, and the isomorphism with $G$ is
one of permutation groups. When $\Fscr$ is such a family, we define:\\
\textbf{Property $\Cbf_{\Fscr,\ell}(\Del)$}
\emph{holds if for all   fields  $K \in \Fscr$, $|\Cl_K[\ell]| \ll_{n,\ell,\Del,\ep} D_K^{\Del+\ep}$ for all $\ep>0$.}\\
Since  Property $\Cbf_{\Fscr,\ell}(\Del)$ remains out of reach for almost all families,  we also consider:\\
\textbf{Property   $\Cbf_{\Fscr,\ell}^*(\Del)$}
\emph{holds if for almost all fields $K \in \Fscr$, $|\Cl_K[\ell]| \ll_{n,\ell,\Del,\ep} D_K^{\Del+\ep}$ for all $\ep>0$.}
We say that a result holds for ``almost all'' fields in a family $\Fscr$ if the subset $E(X)$ of possible exceptions  is density zero   in $\Fscr(X)$, in the sense that 
\[   \frac{|E(X)|}{|\Fscr(X)|}  \maps 0 \qquad \text{as $X \maps \infty$.}\]

Here too, the first progress came for imaginary quadratic fields.
  Soundararajan observed that among  imaginary quadratic fields with discriminant in a dyadic range $[-X,-2X],$ at most one can fail to satisfy $|\Cl_K[\ell]| \ll D_K^{1/2 - 1/2\ell + \ep}$ \cite{Sou00}. This verified $\Cbf_{\Fscr_2^-,\ell}^*(\Del_{\mathrm{GRH}})$   for all primes $\ell$.
 For $\ell=3$ and quadratic fields, Wong observed that $\Cbf_{\Fscr_2^\pm,3}^*(1/4)$ holds \cite{Won01b}.
For any odd prime $\ell$, Heath-Brown and Pierce went below the GRH-bound, proving $\Cbf_{\Fscr_2^-,\ell}^*(1/2 - 3/(2\ell+2))$    \cite{HBP17}. They used the large sieve to show that aside from at most $O(X^\ep)$ exceptions, all discriminants $-d \in [-X,-2X]$ have $|\Cl_{\Q(\sqrt{-d})}[\ell]|$ controlled by counting the number of distinct primes $p_1,p_2$ of a certain size such that (\ref{prime_equation}) has a nontrivial integral solution $(y,z)$. Then they showed there can be few such solutions, while averaging nontrivially over $d$.  These methods relied heavily on the explicit nature of methods for imaginary quadratic fields. Fields of higher degree need a different approach.

\subsection{Dual problems: counting primes, counting fields}

 To apply the Ellenberg--Venkatesh criterion, we face a question such as: ``Given a field, how many small primes split completely in it?''   This question is very difficult in general (and is related to the Generalized Riemann Hypothesis). There is  a  dual question: ``Given a prime, in how many fields does it split completely?'' Ellenberg, Pierce and Wood devised a method to apply the Ellenberg--Venkatesh criterion by tackling the dual question instead \cite{EPW17}. The idea goes like this: suppose   that each   prime splits completely in a positive proportion of fields in a family $\Fscr$. Then the mean number of primes $p \leq x$ that split completely in each field should be comparable to $\pi(x)$, and unless the primes conspire, almost all fields in $\Fscr$ should have close to the mean number of primes split completely in them. To prove that the primes cannot conspire, Ellenberg, Pierce and Wood developed a sieve method, modeled on the Chebyshev inequality from probability. 

The sieve requires as input,   precise counts for the cardinality
\[ N_\Fscr (X;p) = |\{ K \in \Fscr(X): \text{$p$ splits completely in $K$}
     \}|.
\]
It also requires analogous counts $N_\Fscr(X;p,q)$ for when two primes $p \neq q$ split completely in $K$.
Suppose one can prove that for some $\sig>0$ and $\tau < 1$,  for all distinct primes $p,q$,
\beq\label{field_count_pq}
N_\Fscr(X;p,q) = \del(pq) |\Fscr(X)| +O((pq)^\sig |\Fscr(X)|^\tau),
\eeq
for a multiplicative density function $\del(pq)$ taking values in   $(0,1)$.
Then Ellenberg, Pierce, and Wood prove that there exists $\Del_0>0$ (depending on $\tau, \sig$) such that  the mean number of primes $p \leq X^{\Del_0}$ that  split completely in fields in $\Fscr(X)$ is comparable to $\pi (X^{\Del_0})$. Moreover, there can be at most   $O(|\Fscr(X)|^{1- \Del_0})$ exceptional fields $K$ in $\Fscr(X)$ such that fewer than half the mean number of primes split completely in $K$.
 Consequently, for any family $\Fscr$  for which the crucial count (\ref{field_count_pq}) can be proved, combining this sieve with the Ellenberg--Venkatesh criterion proves that   $\Cbf_{\Fscr,\ell}^*(\Del)$ holds for every integer $\ell \geq 2$, where $\Del = \max \{\frac{1}{2}- \Del_0,\Del_{\mathrm{GRH}}\}$.

For which families of fields  can (\ref{field_count_pq}) be proved?
Counting number fields is itself a difficult question. For each integer $D \geq 1$, there are a finite number of extensions $K/\Q$ of degree $n$ and discriminant exactly $D$, by Hermite's finiteness theorem   \cite[\S 4.1]{Ser97}. 
Let $N_n(X)$ denote the number of degree $n$ extensions $K/\Q$   with $D_K \leq X$ (counted up to isomorphism).
A  folk conjecture, sometimes associated to Linnik, states that 
\beq\label{Linnik}
N_n(X) \sim c_n X \qquad \text{as $X \maps \infty$.}
\eeq
When $n=2$ this is essentially equivalent to counting square-free integers (see \cite[Appendix]{EPW17}).
For degree $n=3$, this is a deep result of Davenport and Heilbronn \cite{DavHei71}. For degree $n=4$ it is known by  celebrated results of Cohen, Diaz y Diaz, and Olivier (counting quartic fields $K$ with $\Gal(\tilde{K}/\Q) \simeq D_4$), and Bhargava  (counting non-$D_4$ quartic fields) \cite{CDO02,Bha05}. For degree $n=5$, it is known by landmark work of  Bhargava \cite{Bha10a}.

The sieve method of Ellenberg, Pierce and Wood requires an even more refined count (\ref{field_count_pq}), with prescribed local conditions and a power-saving error term with explicit dependence on $p,q$. 
Power saving error terms for $N_n(X)$ were   found for $n=3$ by Belabas, Bhargava, and Pomerance \cite{BBP10}, Bhargava, Shankar, and Tsimerman \cite{BST13}, Taniguchi and
Thorne \cite{TanTho13}; for $n=4$ (non-$D_4$) by Belabas, Bhargava, and Pomerance \cite{BBP10}; and for $n=5$ by Shankar and Tsimerman \cite{ShaTsi14}.
These results can be refined to prove (\ref{field_count_pq}).
 Ellenberg, Pierce and Wood used this strategy to prove that  when $\Fscr$ is the family of fields of degree $n=2,3,4$ (non-$D_4$), or $5$,  $\Cbf_{\Fscr,\ell}^*(\Del_{\mathrm{GRH}})$ holds  for all sufficiently large primes $\ell$. (For the few remaining small $\ell$, $\Cbf_{\Fscr,\ell}^*(\Del)$ holds with a slightly larger $\Del<1/2$.)  
 Counting  quartic $D_4$-fields with local conditions, ordered by discriminant, remains an interesting open problem.

 The probabilistic method of Ellenberg--Pierce--Wood   uses the property that the density function $\del(pq)$  in (\ref{field_count_pq}) is multiplicative (i.e. local conditions at $p$ and $q$ are asymptotically independent).  
Frei and Widmer have   adapted this approach to prove $\Cbf_{\Fscr,\ell}^*(\Del_{\mathrm{GRH}})$ for all sufficiently large $\ell$, for $\Fscr$ a family of totally ramified cyclic extensions of   $k$ \cite{FreWid18}. (That is, $\Fscr$ is comprised of cyclic extensions $K/k$ of degree $n$ in which every prime ideal of $\Ocal_k$ not dividing $n$ is either unramified or totally ramified in $K$). 
This family is   chosen since the density function $\del(pq)$ is multiplicative.
It would be interesting to investigate whether a probabilistic method can rely less strictly upon multiplicativity of the density function.

There is a great obstacle to expanding the above approach to the family of all fields of degree $n$ when $n \geq 6$. Then, even the   asymptotic (\ref{Linnik}) is not known.
For each $n \geq 6$,  
\beq\label{fields}
N_n(X) \leq a_n X^{c_0 (\log n)^2}
\eeq
is the best-known bound, with $c_0=1.564$, by Lemke Oliver and Thorne \cite{LemTho20}; this improves on Couveignes \cite{Cou20},   Ellenberg and Venkatesh \cite{EllVen06}, and Schmidt \cite{Sch95}.
 For lower bounds, in general the record is    $N_n(X)\gg X^{1/2 + 1/n}$, for all $n \geq 7$ \cite{BSW16x}. For any $n$ divisible by $p=2,3$ or $5$, Kl\"uners (personal communication) has observed that  $N_n(X) \gg X$, since there exists a field $F/\Q$ of degree $n/p$ such that degree $p$  $S_p$-extensions of $F$ exhibit linear asymptotics.

Tackling the problem of counting primes with certain splitting conditions in a specific field via the dual problem of counting fields with certain local conditions at specific primes seems out of reach for higher degree fields. How about tackling the problem of counting primes directly?

\section{Counting primes with $L$-functions}

The prime number theorem states that the number $\pi(x)$  of primes $p \leq x$ satisfies $\pi(x) \sim \Li(x)$ as $x \maps \infty$. To count \emph{small} primes, or primes in short intervals, requires understanding the error term as well as the main term. For each $ 1/2 \leq \Del<1$, the statement
\beq\label{PNT}
\pi(x)  = \Li(x) + O(x^{\Del+\ep}) \qquad \text{for all $\ep>0$} 
\eeq
is essentially equivalent to the statement that the Riemann zeta function $\zeta(s)$   is zero-free for $\Re(s) > \Del$ \cite[Ch. 18]{Dav00}.
%\cite[Ch. 18]{Dav00}
The Riemann Hypothesis conjectures this is true for $\Del=1/2$; it is not known for any $\Del<1$. The best known Vinogradov--Korobov zero-free region is: %$\zeta(\sig +it) \neq 0$ for
\beq\label{VinKor}
\sig \geq 1 - \frac{C}{(\log t)^{2/3} (\log \log t)^{1/3}}, \qquad t \geq 3,
\eeq
  with an absolute constant $C>0$ computed by Ford  \cite{For02}.

To count primes with a specified splitting type in a Galois extension $L/\Q$ of degree $n_L \geq 2$, 
 consider the counting function 
\beq\label{pi_dfn}
 \pi_{\mathscr{C}}(x,L/\Q) = | \{ p \leq x: \text{$p$ unramified in $L$}, \left[ \frac{L/\Q}{p} \right] = \Cscr \}|,
 \eeq
in which $\left[ \frac{L/\Q}{p} \right]$ is the Artin symbol and $\Cscr$ is any fixed conjugacy class in $G=\Gal(L/\Q)$.
For example, when $L=\Q(e^{2\pi i /q})$, this can be used to count primes in a fixed residue class modulo $q$. Or for example, for any Galois extension $L/\Q$, when $\Cscr=\{\mathrm{Id}\}$, this counts primes that split completely in $L$. 
By the celebrated Chebotarev density theorem \cite{Che26},
  \beq\label{Cheb_original}
    \pi_{\mathscr{C}}(x,L/\Q)  \sim \frac{|\Cscr |}{|G|} \mathrm{Li} (x), \qquad \text{as $x \maps \infty$.}
    \eeq
But just as for $\pi(x)$, to count \emph{small} primes accurately requires more quantitative information.
A central goal is to prove an asymptotic for $\pi_{\mathscr{C}}(x,L/\Q)$ that is valid for $x$ very small relative to $D_L = |\mathrm{Disc}L/\Q|$, and with an effective error term.   This requires exhibiting a zero-free region for the Dedekind zeta function $\zeta_L(s)$. 
This is more complicated than (\ref{VinKor}), due to the possibility of an exceptional Landau--Siegel zero: within the region
 \beq\label{LO_region}
 \sig \ge 1-(4\log D_L)^{-1}, \qquad |t| \le (4\log D_L)^{-1},
\eeq
$\zeta_L(\sig + it)$ can contain at most one (real, simple) zero,  denoted $\be_0$ if it exists. (As observed by Heilbronn and generalized by Stark, if $\be_0$ exists then it must ``come from'' a quadratic field, in the sense that $L$ contains a quadratic subfield $F$ with $\zeta_F(\be_0)=0$  \cite{Hei72,Sta74}.)
 
 Lagarias and Odlyzko used the zero-free region (\ref{LO_region}) to prove there exist absolute, computable constants $C_1,C_2$ such that for all  
$x \geq \exp(10n_L(\log D_L)^2),
$
\beq\label{LO_error}
  \left| \pi_{\mathscr{C}}(x,L/\Q)  - \frac{|\mathscr{C}|}{|G|} \Li(x)  \right| \leq  \frac{|\mathscr{C}|}{|G|}  \Li (x^{\be_0}) + C_1 x \exp(-C_2  n_L^{-1/2} (\log x)^{1/2} ),
  \eeq
in which the $\be_0$ term is present only if $\be_0$ exists (see \cite{LagOdl75}, and Serre \cite{Ser82}). This was the first effective Chebotarev density theorem. It can be difficult to apply to  questions of interest because of   the mysterious $\be_0$ term, and because $x$ must be a large power of $D_L$ (certainly at least $x \geq D_L^{10 n_L}$). In contrast,  to apply  the Ellenberg--Venkatesh criterion to a field $K$ of degree $n$, we aim to exhibit primes $p < D_K^{\eta}$ that split completely in the Galois closure $\tilde{K}$ (and hence in $K$), with   $\eta\approx 1/(2\ell(n-1)) \maps 0$ as $n,\ell \maps \infty$.
(These primes are even smaller relative to $D_{\tilde{K}}$, since  $D_K^{|G|/n} \ll_G D_{\tilde{K}} \ll_G D_K^{|G|/2}$, where $G = \Gal(\tilde{K}/\Q)$ \cite{PTBW20}.)

If GRH holds for $\zeta_L(s)$, then $\zeta_L(s)$ is zero-free for $\Re(s)>1/2$, and Lagarias and Odlyzko improve (\ref{LO_error})  in three ways: 
(i) it is valid for $x \geq 2$; (ii) the $\be_0$ term is not present; (iii) the remaining error term is $O(x^{1/2} \log (D_L x^{n_L}))$. Properties (i) and (ii)   show that for every $\eta>0$, for every degree $n$ extension $K/\Q$ with $D_K$ sufficiently large, at least $\gg \pi(D_{K}^\eta)$ primes $p \leq D_{K}^{\eta}$  split completely in the Galois closure $\tilde{K}$ (and hence in $K$). When input into   the Ellenberg--Venkatesh criterion, this is the source of the GRH-bound (\ref{GRH_bound}) for all integers $\ell \geq 2$.

 Here is a central goal: 
 improve the Chebotarev density theorem (\ref{LO_error})   without assuming GRH, so that: (i') for any $\eta>0$ it is valid for $x$ as small as $x \geq D_{L}^\eta$ (for all $D_L$ sufficiently large) and (ii) the $\be_0$ term is not present. (For many applications, the final error term in (\ref{LO_error}) suffices as is.) 
 If this held for $L=\tilde{K}$ the Galois closure of a field $K$, the Ellenberg--Venkatesh criterion would imply the GRH-bound  (\ref{GRH_bound}) for $\ell$-torsion in $\Cl_K$ for all integers $\ell \geq 2$, without assuming GRH. Recently,
 Pierce, Turnage-Butterbaugh and Wood showed  that   the key improvements (i') and (ii) hold if for some $0<\del \leq 1/4$,  $\zeta_{L}(s)/\zeta(s)$ is zero-free for $s=\sig + it$ in the box
\beq\label{box}
1-\del \leq \sig \leq 1, \qquad |t| \leq  \log D_{L}^{2/\del}.
\eeq
Proving this for any particular $L$-function $\zeta_L(s)/\zeta(s)$ of interest is out of reach. 
Instead,  it can be productive to study a family of $L$-functions.  
In particular, if   $\Fscr=\Fscr_n(G;X)$ is a family of degree $n$ fields with fixed Galois group of the Galois closure, property $\Cbf^*_{\Fscr,\ell}(\Del_{\mathrm{GRH}})$ will follow (for all integers $\ell \geq 2$) if it is true   for \emph{almost all} fields $K \in \Fscr_n(G;X)$, that $\zeta_{\tilde{K}}(s)/ \zeta(s)$ is zero-free in the box  (\ref{box}).
This was the strategy  Pierce, Turnage-Butterbaugh and Wood developed in \cite{PTBW20}, which we will now briefly sketch.

\subsection{Families of $L$-functions}

There is a long history of estimating the density of zeroes within a certain region, for a  family of $L$-functions. If we can show there are fewer possible zeroes in the region than there are $L$-functions in the family, then  some of the $L$-functions must be zero-free in that region. We single out a result of Kowalski and Michel, who  used the large sieve to prove a   zero density result for families of cuspidal automorphic $L$-functions \cite{KowMic02}. In particular, for suitable families, their result  implies that almost all $L$-functions in the family must be zero-free in a box analogous to (\ref{box}).

There are two fundamental barriers to applying this to our problem of interest: the representation underlying $\zeta_{\tilde{K}}(s)/\zeta(s)$ is not always cuspidal, and it is not always known to be automorphic. 
Suppose $G$ has irreducible complex representations $\rho_0,  \rho_1, \ldots, \rho_r$, with $\rho_0$  the trivial representation. Then for $K \in \Fscr_n(G;X)$, $\zeta_{\tilde{K}}$ is a product of Artin $L$-functions, 
\beq\label{zeta_prod}
 \zeta_{\tilde{K}}(s)/\zeta(s)  =  \prod_{j=1}^r L(s,\rho_j, \tilde{K}/\Q)^{\mathrm{dim}\rho_j}. \eeq
The Artin (holomorphy) conjecture posits that for each   nontrivial irreducible representation $\rho_j$, $L(s,\rho_j, \tilde{K}/\Q)$ is entire. The (strong) Artin conjecture posits that  for each nontrivial irreducible representation $\rho_j$, there is an associated cuspidal automorphic representation $\pi_{\tilde{K},j}$ of $\GL(m_j)/\Q$, and  $L(s,\pi_{\tilde{K},j}) = L(s,\rho_j,{\tilde{K}}/\Q)$.  This is known for  certain types of representations of certain groups, but otherwise is a deep open problem (see   recent work in \cite{Cal13}). For the moment, we will proceed by assuming the strong conjecture.
Then the factorization (\ref{zeta_prod}) naturally slices the family $\zeta_{\tilde{K}}(s)/\zeta(s)$, as $K$ varies over $\Fscr_n(G;X)$, into $r$ families $\Lscr_1(X), \Lscr_2(X),\ldots, \Lscr_r(X)$, where
  each  $\Lscr_j(X)$ is the set of cuspidal automorphic representations $\pi_{\tilde{K},j}$ associated to the  representation $\rho_j$.   Kowalski and Michel's result applies to each family $\Lscr_j(X)$ individually. This proves that  
 every representation $\pi \in \Lscr_j(X)$ has associated $L$-function $L(s,\pi)$  being zero-free in the box (\ref{box})---except for a possible subset of ``bad'' representations $\pi$, of density zero in $\Lscr_j(X)$, for which $L(s,\pi)$ could have a zero in the box. 
 (Of course, no such zero exists if GRH is true, but we are not assuming GRH.)

Now a crucial difficulty arises: if there were a ``bad'' representation $\pi \in \Lscr_j(X)$, in how many products (\ref{zeta_prod}) could it appear, as $K$ varies over $\Fscr_{n}(G;X)$? Each field $K$ for which the ``bad'' factor $L(s,\pi)$ appears could have a zero of $\zeta_{\tilde{K}}(s)/\zeta(s)$ in (\ref{box}).
Thus the crucial question is:  for a fixed nontrivial irreducible representation $\rho$ of $G$,  how many   fields $K_1,K_2  \in \Fscr_n(G;X)$   have 
$
 L(s,\rho, \tilde{K}_1/\Q)  = L(s,\rho, \tilde{K}_2/\Q)?
$
This  can be stated a different way. Given a subgroup $H$ of $G$, let $\tilde{K}^H$ denote the subfield of $\tilde{K}$ fixed by $H$. 
It turns out that the question can be transformed into:
 how many fields $K_1,K_2 \in \Fscr_n(G;X)$ have  $\tilde{K}_1^{\Ker(\rho)} = \tilde{K}_2^{\Ker(\rho)}$? Let us call this a collision. If a positive proportion of fields in $\Fscr_n(G;X)$ can collide for $\rho_j$, then via the factorization (\ref{zeta_prod}), the possible existence of even one ``bad'' element in $\Lscr_j(X)$ could allow a positive proportion of the functions $\zeta_{\tilde{K}}(s)/\zeta(s)$ to have a zero in  (\ref{box}). In particular, then this approach would fail to prove $\Cbf^*_{\Fscr,\ell}(\Del_{\mathrm{GRH}})$ for the family $\Fscr = \Fscr_n(G;X).$
 To rule this out, we   aim to show that for each nontrivial irreducible representation $\rho_j$ of $G$, collisions are rare.
 
We define the ``collision problem'' for the family $\Fscr_n(G;X)$: how big is
\beq\label{collisions}
 \max_{\rho} \max_{K_1 \in \Fscr_n(G;X)}|\{ K_2 \in \Fscr_n(G;X): \tilde{K_1}^{\Ker(\rho)} = \tilde{K_2}^{\Ker(\rho)}\}|?
\eeq
Here the maximum is over the nontrivial irreducible representations $\rho$ of $G$ with $\Ker(\rho)$   a proper normal subgroup of $G$.
Suppose for a particular family $\Fscr_n(G;X)$, the collisions (\ref{collisions}) number at most $\ll X^\al$.  Then the strategy sketched here  ultimately shows that aside from at most $\ll X^{\al+\ep}$ exceptional fields (for any $\ep>0$), every field in $K \in \Fscr_n(G;X)$ has the property that  an improved Chebotarev density theorem with properties (i') and (ii) holds for its Galois closure $\tilde{K}$. If we can prove simultaneously that $|\Fscr_n(G;X)| \gg X^\be$ for some $\be>\al$, then the improved Chebotarev density theorem holds for almost all fields in the family.
Consequently, we would obtain property $\Cbf^*_{\Fscr,\ell}(\Del_{\mathrm{GRH}})$, for all integers $\ell \geq 2$.

Thus the goal of bounding $\ell$-torsion in class groups of fields in the family $\Fscr_n(G;X)$ has been transformed into a question of counting how often certain fields share a subfield.
For which families can the collision problem (\ref{collisions}) be controlled?
For some groups, the number of collisions can be $\gg|\Fscr_n(G;X)|$ (for example $G = \Z/4\Z)$. On the other hand, if $G$ is a simple group, or if all nontrivial irreducible representations of $G$ are faithful,  the number of collisions is $\ll1$ (but a lower bound $|\Fscr_n(G;X)|\gg X^\be$ for some $\be>0$   may not yet be known). In general, controlling the collision problem is difficult. 

One idea is to restrict attention to an advantageously chosen subfamily of fields, call it $\Fscr_n^*(G;X) \subset \Fscr_n(G;X)$. To bound (\ref{collisions}) within a subfamily   it  suffices to count 
\beq\label{restricted_collisions}
\max_H \max_{\substack{F\\ \deg(F/\Q)=[G:H]}} |\{ K \in \Fscr_n^*(G;X) :\tilde{K}^{H}=F  \}| .
\eeq
Here $H$ ranges over the proper normal subgroups of $G$ that appear as the kernel of some nontrivial irreducible representation. For some groups $G$, if $\Fscr_n^*(G;X)$ is defined appropriately, this   can be further transformed into counting number fields with $\emph{fixed}$ discriminant.

Let us see how this goes in the example  $G = S_n$ with $n=3$ or $n \geq 5$, so that $A_n$ is the only nontrivial proper normal subgroup  (the kernel of the sign representation). 
Consider the subfamily $\Fscr_n^*(S_n;X)$ of fields with square-free discriminant. (These are a positive proportion of all degree $n$ $S_n$-fields for $n \leq 5$ and conjecturally so for $n \geq 6$.) %\cite{DavHei71,Mal02,Bha05,Bha07,Bha10a,Bha14x} 
Then for $H=A_n$ and $F$ a fixed quadratic field, it can be shown that any field $K$ counted in (\ref{restricted_collisions}) must have the property that $D_K=D_F$ (up to some easily controlled behavior of wildly ramified primes). Under this very   strong identity of discriminants,   (\ref{restricted_collisions}) is dominated by 
\beq\label{collision_fixed}
\max_{D \geq 1} |\{K \in \Fscr_n^*(S_n;X) : D_K =D\}|.
\eeq
This strategy transforms the collision problem into counting   fields of \emph{fixed} discriminant.  

For certain other   groups $G$, (\ref{restricted_collisions}) can also be dominated by a quantity analogous to (\ref{collision_fixed}) if  the subfamily $\Fscr_n^*(G;X)$ is defined by specifying that each prime that is tamely ramified in $K$ has its inertia group generated by an element in a carefully chosen conjugacy class $\Iscr$ of $G$.
For such a group $G$, the final step in this strategy for  proving improved Chebotarev density theorems for almost all fields in the family $\Fscr_n^*(G;X)$ is to bound (\ref{collision_fixed}). If $|\Fscr_n^*(G;X)| \gg X^\be$ is known, it suffices to prove (\ref{collision_fixed}) is $\ll X^\al$ for some $\al< \be$. In general, counting number fields with fixed discriminant   is very difficult---we will return to this problem later. But for some families $\Fscr_n^*(G;X)$, (\ref{collision_fixed}) can be controlled sufficiently well, relative to a known lower bound for $|\Fscr_n^*(G;X)|.$

This is the strategy developed by Pierce, Turnage-Butterbaugh and Wood in \cite{PTBW20}.
 The result is an improved   Chebotarev density theorem, with properties (i') and (ii), that holds unconditionally for almost all fields in the following families:
(a) $\Fscr_p(C_p;X)$ cyclic extensions of any prime degree; (b) $\Fscr_n^*(C_n;X)$ totally ramified cyclic extensions of any degree $n \geq 2$;
(c) $\Fscr_p^*(D_p;X)$ prime degree dihedral extensions, $\Iscr$ being the class of order 2 elements; 
(d) $\Fscr_n^*(S_n;X)$ fields of square-free discriminant, $n=3,4$;
(e) $\Fscr_4^*(A_4;X)$, $\Iscr$ being either class of order 3 elements. 
 Conditional on the strong Artin conjecture, they proved the improved Chebotarev density theorem also holds for almost all fields in the following families:
(f) $\Fscr_5^*(S_5;X)$ quintic fields of square-free discriminant;
(g) $\Fscr_n(A_n;X)$, for all $n \geq 5.$
(There are other families, such as $\Fscr_n^*(S_n;X)$ for $n \geq 6$, to which the strategy   applies, but the current upper bound known for (\ref{collision_fixed}) is larger than the known lower bound  for $|\Fscr_n^*(S_n;X)|$.)
As a consequence, Pierce, Turnage-Butterbaugh and Wood proved for each family (a)-(e) that $\Cbf_{\Fscr,n}^*(\Del_{\mathrm{GRH}})$ holds unconditionally for all integers $\ell \geq 2$,  and it holds for each family (f)-(g) under the strong Artin conjecture. This was the first time  such a result was proved for families of fields of arbitrarily large degree.

\subsection{Further developments}
Since the work outlined above, many interesting   new developments have followed, relating to   zero density results  for families of $L$-functions, Chebotarev density theorems for families of fields, and $\ell$-torsion in class groups of fields in specific families.

First, there has been renewed interest in
    zero density results for families of $L$-functions, concerning potential zeroes   in regions close to the line $\Re(s)=1$, and extending the perspective of Kowalski and Michel \cite{KowMic02}; see for example  \cite{ThoZam21a,BTZ19x,HumTho21x}.

Second, several new strategies have focused  on the problem of proving effective Chebotarev density theorems for almost all fields in a family. The work in \cite{PTBW20} raised several desiderata. Some groups $G$ have the property that no ramification restriction exists that allows the ``collision problem'' in the form (\ref{restricted_collisions}) to be transformed into a ``discriminant multiplicity problem'' in the form (\ref{collision_fixed}).  For example, this occurs for any non-cyclic abelian group, or $D_4$. These cases remain open; instead An recently proved a Chebotarev density theorem for almost all fields in a family  of quartic $D_4$-fields associated to a fixed biquadratic field  \cite{An20}. Another significant desideratum was to remove the dependence on the strong Artin conjecture.
Thorner and Zaman recently achieved this, by proving a zero density estimate directly for Dedekind zeta functions, without passing through the factorization (\ref{zeta_prod}) \cite{ThoZam19x}.   But that work is still explicitly conditional on the ability to control a collision problem similar to (\ref{collisions}), for which the best known strategy is still the approach of \cite{PTBW20}.

Most recently, the collision problem   has been bypassed for certain groups $G$ by interesting new work of Lemke Oliver, Thorner and Zaman \cite{LOTZ21x}.
 Their key idea when studying fields in a family $\Fscr_n(G;X)$ is to prove a zero-free region not for $\zeta_{\tilde{K}}/\zeta$ but  for $\zeta_{\tilde{K}}/\zeta_{\tilde{K}^N}$ where  $N$ is a nontrivial normal subgroup of $G$. This allows them to replace a collision problem like (\ref{collisions})  by an ``intersection multiplicity problem,'' bounding
\beq\label{collissions_LOTZ_N}
\max_{K_1 \in \Fscr_n(G;X)}
    |\{ K_2 \in \Fscr_n(G;X): \tilde{K}_1 \intersect \tilde{K}_2 \neq \tilde{K}_1^N \intersect \tilde{K}_2^N\}|.
\eeq
The number of exceptional fields, for which a desired Chebotarev-type theorem cannot be verified, is then dominated  by (\ref{collissions_LOTZ_N}) (up to $X^\ep$).
This is advantageous if $G$ has a unique minimal nontrivial normal subgroup $N$, so that (\ref{collissions_LOTZ_N}) is $\ll 1$.
 But as a trade-off, one no longer obtains an effective Chebotarev density theorem for each conjugacy class $\Cscr$ in $G$.

 Let $\pi_K(x)$ count prime ideals $\pfr \subset \Ocal_K$ with $\Nfr_{K/\Q} \pfr \leq x$. 
  Let $\Fscr$ represent either of the two following families: degree $p$ fields $K/\Q$ for $p$ prime, or degree $n$ $S_n$-fields $K/\Q$, for any $n \geq 2$.  Lemke Oliver, Thorner and Zaman prove that except for at most $\ll X^\ep$ exceptional fields, every $K \in \Fscr(X)$ has $|\pi_K(x) -\pi(x)| \leq C_1 x \exp(-C_2 \sqrt{\log x})$ for every $x \geq (\log D_K)^{C_3(n,\ep)}$. In either family $\Fscr$, they obtain results on $\ell$-torsion    by applying the Ellenberg--Venkatesh criterion using prime ideals of degree 1.
If $\pi_K^*(x)$ counts only prime ideals of degree 1, then
$\pi_K^*(x) =\pi_K(x)+O_n(\sqrt{x})$, so the above result exhibits    many  small prime ideals of degree 1. Thus for either family, $\Cbf_{\Fscr,n}^*(\Del_{\mathrm{GRH}})$
 holds unconditionally for all $\ell$ (and the exceptional set is very small). 
(They also exhibit infinitely many degree $n$ $S_n$-fields $K$ with $\Cl_K$ as large as possible, but $|\Cl_K[\ell]|$ bounded  by (\ref{GRH_bound}) for all $\ell$; and infinitely many totally real degree $n$ $S_n$-fields $K$ with $\Cl_K$ containing an element of exact order $\ell$ and $|\Cl_K[\ell]|$ bounded  by (\ref{GRH_bound}).)
What happens when $G$ does not have a unique minimal nontrivial normal subgroup?
Here is an open question: in general, when $N$ is a nontrivial normal subgroup of $G$ (not necessarily unique or minimal), what is the true order of growth of (\ref{collissions_LOTZ_N}) as $X \maps \infty$?   Questions about this ``intersection multiplicity'' are gathered in \cite{LOTZ21x}.

 Third, increased attention has turned to bounding $\ell$-torsion in class groups for \emph{all} fields in special families specified by the Galois group: that is, proving property $\Cbf_{\Fscr,\ell}(\Del)$ for some $\Del<1/2$. First,  Kl\"{u}ners and Wang have proved   $\Cbf_{\Fscr,p}(0)$ for the family $\Fscr_{p^r}(G;X)$ for any $p$-group $G$;
 this generalizes the application of genus theory to prove $\Cbf_{2,2}(0)$  \cite{KluWan20x}.
 Second, 
 let $G = (\Z/p\Z)^r$ be an elementary abelian group of rank $r \geq 2$, with $p$ prime.   Wang has shown that for every $\ell$, within the family   of Galois $G$-fields $K/\Q$, property $\Cbf_{\Fscr,\ell}(1/2 - \del(\ell,p))$ holds for some $\del(\ell,p)>0$ \cite{Wan21}. Since the savings $\del(\ell,p)$ is independent of the rank, for $r$ sufficiently large    this is better than $\Cbf_{\Fscr,\ell}(\Del_{\mathrm{GRH}})$. The method of proof plays off the interaction of three facts arising from the precise structure of $G$: first, $|\Cl_K[\ell]|$ factors as a product of $|\Cl_{F}[\ell]|$ where $F$ varies over the $\approx p^{r-1}$ many degree $p$ subfields of $K$, so it suffices to bound one of these factors nontrivially. Second,  
 any rational prime splits completely in $\approx p^{r-2}$ of these subfields, so at least one subfield has a positive proportion of   primes splitting completely in it.
 Third, the sizes of the  discriminants of the subfields can be played against each other, so that known prime-counting results (which may \emph{a priori} seem to count primes that are ``too large'') suffice for the application of the Ellenberg--Venkatesh criterion. This is an interesting counterpoint to the methods described earlier.
In another direction, Wang has developed the notion of a forcing extension; certain nilpotent groups can be built from elementary $p$-groups via forcing extensions.  If $G'$ is constructed from $G$ by a forcing extension, then   $\Cbf_{\Fscr',\ell}(\Del')$ can be deduced from $\Cbf_{\Fscr,\ell}(\Del),$ for some $\Del,\Del'<1/2$, where $\Fscr$ is the family of   $G$-extensions and $\Fscr'$ is the family of   $G'$-extensions \cite{Wan20x}.

All of the results mentioned in this section (except where genus theory suffices)  directly apply or build on the Ellenberg--Venkatesh criterion. Can this criterion be strengthened?  Ellenberg has suggested some possible improvements in \cite{Ell08}. In particular, let $
\eta(K):=\inf\{ H_K(\al): K = \Q(\al)\}  
$
denote 
the minimum (relative) multiplicative Weil height of a generating element of $K$.
Roughly speaking, Ellenberg notes the criterion (\ref{EV_criterion}) can actually allow prime ideals with norms as large as $\eta(K)^{1/\ell}$.
The restriction  to norms $<D_K^{\frac{1}{2\ell(n-1)}}$ in (\ref{EV_criterion}) was made since the lower bound $\eta(K) \geq D_K^{\frac{1}{2(n-1)}}$ holds for all fields \cite{Sil84}.
Widmer, also with Frei, has shown that $\eta(K)$ can be enlarged for almost all fields in certain families, leading to improved upper bounds for $\ell$-torsion in those fields \cite{Wid18,FreWid21}.  That is, they improve the very notion of the ``GRH-bound'' (\ref{GRH_bound}), and show that the parameter we have called $\Del_{\mathrm{GRH}}$ can actually be taken smaller for some fields.
Their work raises interesting open questions:  what upper and lower bounds hold for $\eta(K)$, for all (or almost all) fields in a family? Ruppert \cite{Rup98} has conjectured uniform upper bounds $\eta(K) \leq D_K^{1/2}$ (now proved for almost all fields in some families by \cite{PTBW20}). If this is true, the Ellenberg--Venkatesh criterion would hit a barrier, for most fields, with a result like $|\Cl_K[\ell]| \ll D_K^{1/2 -1/2\ell+\ep}$ for any degree $n$, still far from the $\ell$-torsion Conjecture.
It would be very  interesting to find a new, different criterion.

 \section{Why do we expect the $\ell$-torsion Conjecture   to be true?}
Recall that the $\ell$-torsion Conjecture \ref{conj_class}  is still known only in the case stemming from Gauss's work: $n=2, \ell=2$. It is a good idea to affirm why we   believe the $\ell$-torsion Conjecture should be true. We will consider this from three perspectives.

\subsection{From the perspective of the Cohen--Lenstra--Martinet heuristics}

So far, when we have mentioned a result  for almost all fields in a family,  we have not focused   on the size of a potential exceptional set, other than showing it is smaller than the size of the full family. But to understand the $\ell$-torsion Conjecture, we must quantify a potential exceptional set, and show that for all sufficiently large discriminants, it is empty.

Let us abstract this, for a family $\Fscr_0(X)$ of fields $K$ with $D_K$ in a dyadic range $(X/2,X]$, from which more general results can easily be deduced by summing over $\ll \log X$ dyadic ranges.
Suppose $f: \Fscr_0(X) \maps \N$ is a function with   $f(K) \leq D_K^a$ for all $K$. Suppose that for some $\Del<a$ we can improve this to $f(K) \leq D_K^\Del$ for all $K$ outside of some exceptional set $ E_0^\Del(X) \subset \Fscr_0(X).$ Then 
 \beq\label{average_sum}  \sum_{K \in \Fscr_0(X) } f(K)= \sum_{K \in \Fscr_0(X) \setminus E_0^\Del(X)} f(K) + \sum_{K \in E_0^\Del(X)} f(K)
 \leq |\Fscr_0(X)|X^{\Del} + |E_0^\Del(X)| X^a. 
 \eeq
As long as  $|E_0^\Del(X)| \ll |\Fscr_0(X)|X^{ - (a-\Del)},$ this shows   that $f(K) \ll X^\Del$ on average.
On the other hand, suppose we know $\sum_{K \in \Fscr_0(X)} f(K) \leq X^b$. Then a potential set of exceptions $E_0^\Del(X)= \{ K \in \Fscr_0(X): f(K) > D_K^\Del\}$ can be controlled by
\beq\label{moment}
X^\Del |E_0^\Del(X)| \ll \sum_{K \in E_0^\Del(X)} f(K) \leq \sum_{ K \in \Fscr_0(X)  } f(K) \leq  X^b.
\eeq
Thus   $|E_0^\Del(X)| \ll X^{b - \Del}$, and exceptional fields are density zero in $\Fscr_0(X)$,  provided $X^{b-\Del} = o(|\Fscr_0(X)|)$.
That is,  a nontrivial upper bound on $\ell$-torsion  for ``almost all'' fields in a family $\Fscr$ is essentially equivalent to   the same upper bound  ``on average.''

To verify the $\ell$-torsion Conjecture we wish to show a ``pointwise'' bound:  for every $\ep>0$, there exists $D_\ep$ such that when $D_K \geq D_\ep,$ there are \emph{no} exceptions to the bound $|\Cl_K[\ell]| \leq D_K^\ep.$
The key is to consider not averages but arbitrarily high $k$-th moments.
In the general setting above, suppose that we know   $\sum_{K \in \Fscr_0(X)} f(K)^k \leq X^b$, for a real number $k \geq 1$. Then for any fixed $\Del>0$, adapting the argument (\ref{moment}) shows that 
$|E_0^\Del(X)| \ll X^{b-k\Del}$. If  the $k$-th moment is uniformly bounded by $X^b$ for a sequence of $k \maps \infty$, then for each $\Del>0$, we can take $k$ sufficiently large to conclude that the set of exceptions is empty. 

This perspective has been applied by Pierce, Turnage-Butterbaugh and Wood in \cite{PTBW21} to prove that the $\ell$-torsion Conjecture holds for \emph{all} fields in a family $\Fscr(X)$
if there is a real number $\al \geq 1$ such that for a sequence of arbitrarily large $k$,
  \beq\label{moment_bound}
 \sum_{K \in  \Fscr(X)}  |\Cl_K[\ell]|^k \ll _{n,\ell, k,\al} |\Fscr (X)|^{\al}, \qquad \text{for all $X \geq 1$.}
 \eeq
The Cohen--Lenstra--Martinet heuristics predict that (\ref{moment_bound}) holds, in the form of an even  stronger asymptotic with $\al=1$, for all integers $k \geq 1$, for families of Galois $G$-extensions, at least for all primes $\ell \ndiv |G|$.   The appropriate moment formulation can be found in \cite{CohLen84} for degree 2 fields and in \cite{WanWoo21} for higher degrees, building on \cite{CohMar90}.
This confirms that the $\ell$-torsion Conjecture follows from another well-known set of conjectures. 

The Cohen--Lenstra--Martinet heuristics are   a subject of intense interest and much recent activity. Here are some spectacular successes most closely related to our topic.
  Davenport and Heilbronn \cite{DavHei71} have proved
\beq\label{DavHei23} \sum_{\bstack{\deg(K)=2}{0<D_K\leq X}} |\Cl_{K}[3]| \sim \left( \frac{2}{3\zeta(2)} + \frac{1}{\zeta(2)}\right) X;\eeq
 second order terms have been found in \cite{BBP10,BST13,TanTho13}.
  Bhargava  \cite{Bha05}  has proved 
\beq\label{Bha32} \sum_{\bstack{\deg(K)=3}{0<D_K\leq X}}| \Cl_{K}[2] |\sim \left( \frac{5}{48\zeta(3)} + \frac{3}{8\zeta(3)} \right) X,
\eeq
in which each isomorphism class of fields is counted once. Very recently, \cite{LOWW21x} obtained analogues of (\ref{DavHei23}) for averages over $\Fscr_{2^m}(G;X)$ for any permutation group $G \subset S_{2^m}$ that is a transitive permutation 2-group containing a transposition.
See also the work of
Smith on the distribution of  $2^{k}$-class groups in imaginary quadratic fields \cite{Smi17x}; 
Koymans and Pagano on $\ell^k$-class groups of degree $\ell$ cyclic fields \cite{KoyPag18x}; Klys on moments of $p$-torsion in cyclic degree $p$ fields (conditional on GRH for $p \geq 5$)  \cite{Kly20}; Milovic and Koymans on 16-rank in  quadratic fields \cite{KoyMil19,KoyMil21};
 Bhargava and Varma \cite{BhaVar15,BhaVar16} elaborating on (\ref{DavHei23}) and (\ref{Bha32}).

The perspective of moments in (\ref{moment_bound}) provides a strong motivation to prove $k$-th moment bounds for $\ell$-torsion.
Fouvry and Kl\"{u}ners  have proved an asymptotic for  $k$-th moments related to  $4$-torsion when $K$ is quadratic, for all integers $k \geq 1$  \cite{FouKlu06}. 
Heath-Brown and Pierce have proved nontrivial bounds for $k$-th moments of $\ell$-torsion for imaginary quadratic fields, for all odd primes $\ell$ \cite{HBP17}. For example, they establish second moment bounds
\beq\label{quad_moments}
\sum_{\substack{K = \Q(\sqrt{\pm D})\\D \leq X}}
    |\Cl_K[3]|^2 \ll X^{23/18}, 
\qquad
\sum_{\substack{K = \Q(\sqrt{-D}) \\D \leq X}}
    |\Cl_K[\ell]|^2 \ll X^{2 - \frac{3}{\ell+1}}, \quad \text{$\ell \geq 5$ prime,}
\eeq
as well as  results for  $k$-th moments for all $k \geq 1$. 
In general, proving tighter control on the size of an exceptional family $E_0^\Delta(X)$ can be used to deduce a better moment bound for $|\Cl_K[\ell]|$, similar to (\ref{average_sum}). This has recently been exploited by Frei and Widmer, in combination with     refinements of the  Ellenberg--Venkatesh criterion, to improve moment bounds on $\ell$-torsion for the  families of fields   studied in \cite{PTBW20} (if $\ell$ is sufficiently large); see \cite{FreWid21}.

Let us mention a  connection to elliptic curves; this was after all the setting in which Brumer and Silverman initially posed the $\ell$-torsion Conjecture. 
 Let $E(q)$ denote the number of isomorphism classes of elliptic curves over $\Q$ with conductor $q$.
 Brumer and Silverman have  conjectured that 
$E(q) \ll_\ep q^{\ep}$ for every $q \geq 1$,   $\ep>0$ \cite{BruSil96}.
Conditionally, this follows from GRH combined with a weak form of the Birch--Swinnerton-Dyer conjecture. They also showed this follows from the $3$-torsion Conjecture for quadratic fields, by proving
\beq\label{EqH3}
E(q) \ll_\ep q^\ep \max_{1 \leq D \leq 1728q} |\Cl_{\Q(\sqrt{\pm D})}[3]|, \qquad \text{for all $\ep>0$.}
\eeq
 Duke and Kowalski have combined this with the celebrated asymptotic (\ref{DavHei23}) to bound $\sum_{1 \leq q \leq Q} E(q)\ll Q^{1+\ep}$ for every $\ep>0$ \cite{DukKow00}. (See also \cite{FNT92} for ordering by discriminant.) Pierce, Turnage-Butterbaugh and Wood have recently  proved that for all $k \geq 1$, the $k$-th moment of $3$-torsion in quadratic fields  dominates the $\ga k$-th moment of $E(q)$, for a numerical constant $\ga \approx 1.9745...$ coming from \cite{HelVen06}, which sharpened the relation (\ref{EqH3}).
 Thus new moment bounds for $E(q)$ can be obtained from (\ref{quad_moments}), for example.
Here is an open problem: prove that 
$\sum_{1 \leq q \leq Q} E(q)=o(Q).$ This would show for the first time that  integers that are the conductor of an elliptic curve have  density zero in $\Z$. In fact, it is conjectured by Watkins that this average is asymptotic to $cQ^{5/6}$ for a certain constant $c$ \cite{Wat08} (building on an analogous conjecture by Brumer--McGuinness for ordering by discriminant \cite{BruMcG90}). 

To conclude, in this section we saw that the truth of the $\ell$-torsion Conjecture is implied by the truth of the well-known Cohen--Lenstra--Martinet heuristics on the distribution of class groups.

\subsection{From the perspective of counting number fields of fixed discriminant}

Let $K/\Q$ be a degree $n$ extension. The Hilbert class field $H_K$ is the maximal abelian unramified   extension of $K$, and $\Cl_K$ is isomorphic to  $\Gal(H_K/K)$.  A second way to motivate the $\ell$-torsion Conjecture is to count intermediate fields   between $K$ and $H_K.$

Here is an argument recorded by  Pierce, Turnage-Butterbaugh and Wood in \cite{PTBW21}. 
Fix a prime  $\ell$ and write $\Cl_K$ additively, so that  $\Cl_K [\ell] \simeq \Cl_K / \ell \Cl_K.$ Now define the fixed field $L=H_K^{\ell \Cl_K}$ lying between $K$ and $H_K$, so $\Gal(L/K) \simeq \Cl_K[\ell]$. Each surjection  $\Cl_K[\ell] \maps \Z/\ell\Z $ generates an intermediate field $M$, with $K \subset M \subset L$ and $\deg(M/\Q)=n\ell$. If $|\Cl_K[\ell]|= \ell^r$,  say, this produces $\approx \ell^{r-1}$  such fields $M$. The crucial point is that since $H_K$ is an unramified extension,  all these fields satisfy a rigid discriminant identity: $D_M = D_K^\ell.$
Consequently, if we can count how many number fields of degree $n\ell$ can share the same \emph{fixed} discriminant, then we can  bound  $\ell$-torsion in $\Cl_K.$  (We have seen this problem before.)  We   formalize the problem of counting number fields of fixed discriminant as follows:\\
\textbf{Property $\Dbf_n(\Del)$:}
\emph{Fix a degree $ n \geq 2$.  Property ${\Dbf}_n(\Delta)$ holds if  for every $\ep>0$ 
and for every fixed integer $D>1$, at most $\ll_{n,\ep} D^{\Del+\ep}$ fields $K/\Q$ of degree $n$ have $D_K=D$.}\\
The strategy sketched above ultimately proves that property $\Dbf_{n\ell}(\Del)$ implies     $\Cbf_{n,\ell}(\ell\Del)$. 
This leads inevitably to the question: is property $\Dbf_{n\ell}(0)$ true? Here is a conjecture:

 \begin{conjecture}[Discriminant Multiplicity Conjecture]\label{conj_DMC}
For each $n \geq 2$, for every $\ep>0$  
and for every integer $D>1$, at most $\ll_{n,\ep} D^{\ep}$ fields $K/\Q$ of degree $n$ have $D_K=D$.
\end{conjecture}
 This conjecture has been recorded by Duke \cite{Duk98}. It implies the $\ell$-torsion Conjecture, a link noted in \cite{Duk98,EllVen05} and quantified in \cite{PTBW21}. 
 Recall the conjecture (\ref{Linnik}) for counting all fields of degree $n$ and discriminant $D_K \leq X$.
 The Discriminant Multiplicity Conjecture  for degree $n$ would   immediately imply $N_n(X)\ll X^{1+\ep}$, which indicates its level of difficulty. Of course, in general, property $\Dbf_{n}(\Del)$ implies $N_n(X) \ll X^{1+\Del+\ep}$ for all $\ep>0$.
(In terms of lower bounds, Ellenberg and Venkatesh  have noted   there can be $\gg D^{c/ \log \log D}$  extensions $K/\Q$ with a fixed Galois group and fixed discriminant $D$ \cite{EllVen05}.)

  The Discriminant Multiplicity Conjecture posits that $\Dbf_n(0)$ holds for each $n \geq 2$. This is true for $n=2$, but it is not known for any other degree. For 
degrees $n=3,4,5$ the best-known results currently are  $\Dbf_3(1/3)$ by \cite{EllVen07};
$\Dbf_4(1/2)$   as found in \cite{Wri89,Klu06b,PTBW20,PTBW21};  $\Dbf_5(199/200)$  as found in \cite{EPW17}, building on \cite{Bha10a,ShaTsi14}. Currently for $n \geq 6,$
the only result for $\Dbf_n(\Del)$ is a trivial consequence of counting fields of bounded discriminant, as in (\ref{fields}), so in particular $\Del=c_0 (\log n)^2>1$ in those cases. 
It would be very interesting to improve the exponent known for $\Dbf_n(\Del)$, for any fixed degree $n \geq 3$. 

As is the case for many of the problems surveyed in this paper, it can also be profitable to study the problem within a family $\Fscr$ of degree $n$ extensions:\\
\textbf{Property $\Dbf_{\Fscr,n}(\Del)$:}
\emph{Fix a degree $n \geq 2$.  Property ${\Dbf}_{\Fscr,n}(\Delta)$ holds if  for every $\ep>0$ 
and for every fixed integer $D>1$, at most $\ll_{n,\ep} D^{\Del+\ep}$ fields $K/\Q$ in the family $\Fscr$ have $D_K=D$.}
This is the type of property Pierce, Turnage-Butterbaugh and Wood used to control the collision problem, in the form (\ref{collision_fixed})  \cite{PTBW20}.
Property $\Dbf_{\Fscr,n}(0)$ has recently been proved by Kl\"uners and Wang, for the family $\Fscr=\Fscr_n(G;X)$ of degree $n$ $G$-extensions for any nilpotent group $G$. This was  built   from the truth of property $\Cbf_{\Fscr,p}(0)$ for $\Fscr$ being the family of Galois $H$-extensions for $H$ a $p$-group, in \cite{KluWan20x}.
There are many other cases where it is an interesting open problem to improve the known bound for Property $\Dbf_{\Fscr,n}(\Del)$.

To conclude, in this section we saw that the $\ell$-torsion Conjecture follows from the Discriminant Multiplicity Conjecture. Now, recall that we saw in the context of bounding $\ell$-torsion that uniform bounds for arbitrarily high moments can imply strong ``pointwise'' results for \emph{every} field.
Can the method of moments be used to approach the Discriminant Multiplicity Conjecture too? 
We turn to this idea next.

\subsection{From the perspective of counting number fields of bounded discriminant}
We come to a  third motivation to believe the $\ell$-torsion Conjecture.  
 Recall the definition (\ref{family_dfn}) of a family  $\Fscr_n(G;X)$ of degree $n$ fields $K/\Q$ with $\Gal(\tilde{K}/\Q)$   isomorphic (as a permutation group) to a nontrivial transitive subgroup $G \subseteq S_n$.   Each element $g \in G$ has an  index defined by
$\mathrm{ind}(g) = n - o_g$, where $o_g$ is the number of orbits of $g$ when it acts on a set of $n$ elements. Define $a(G)$ according to $a(G)^{-1} = \min \{ \mathrm{ind}(g) : 1 \neq g \in G\}$; we see that $\frac{1}{n-1}\leq a(G) \leq 1$. 
Malle has made a well-known conjecture \cite{Mal02}:
\begin{conjecture}[Malle]
For
each $n \geq 2$, for each transitive subgroup $G \subseteq S_n$,
\beq\label{weak_Malle}
 |\Fscr_n(G;X)| \ll_{G,\ep} X^{a(G) + \ep}, \qquad \text{for all $\ep>0$.}
\eeq
Also, $|\Fscr_n(G;X)| \gg_G X^{a(G)}$.
\end{conjecture} 
The full statement of this conjecture is an open problem. Its difficulty is indicated by the fact that it implies a positive solution to the inverse Galois problem for number fields. (A refinement in \cite{Mal04}   specified a power of $\log X$ in place of $X^\ep$; counterexamples to this refinement have been found in
 \cite{Klu05}, but the upper bound  in (\ref{weak_Malle}) is expected to be true.)

 Malle's Conjecture has been proved  for abelian groups, with a strategy by Cohn \cite{Coh54}, and asymptotic counts by M\"{a}ki \cite{Mak85}, Wright \cite{Wri89}.
 For $n=3,4,5$ it is known for $S_n$ by the asymptotic (\ref{Linnik}),
and for $D_4$ by   
Baily \cite{Bai80} (refined to an asymptotic in \cite{CDO02}).
It is known for  $C_2 \wr H$ under mild conditions on $H$ (in particular for at least one group of order $n$ for every even $n$) by \cite{Klu12}, and for $S_n \times A$ with $A$ an abelian group ($n=3,4,5$) by \cite{Wan21a,MTTW20x}.
For prime degree $p$ $D_p$-fields, upper and lower bounds are closely related to $p$-torsion in class groups of quadratic fields, and have been studied in \cite{Klu06,CohTho20,FreWid21}.

For many  groups, it is a difficult open problem to prove upper or lower bounds approaching Malle's prediction.
In many results surveyed here, proving a lower bound for $|\Fscr_n(G;X)|$ has been an important step,  to verify a result applies to ``almost all'' fields in a family. 
For many groups $G$, it is not even known that $|\Fscr_n(G;X)|\gg X^\be$ for some $\be>0$ as $X \maps \infty.$
Here is a tool to prove such a result: 
suppose $f(X,T_1,\dots,T_s) \in \Q[X,T_1,\ldots, X_s]$ is a regular polynomial of total degree $d$ in the $T_i$ and of degree $m$ in $X$ with transitive Galois group $G\subset S_n$ over $\Q(T_1,\dots,T_s)$.  
Then $|\Fscr_n(G;X)| \gg_{f,\ep}X^{\be-\ep}$ for every $\ep>0$, with $\be =\frac{1-|G|^{-1}}{d(2m-2)}$; this is proved in \cite{PTBW20}.
For $G=A_n$, a polynomial $f$ exhibited by Hilbert can be input to this criterion, implying that $|\Fscr_n(A_n;X)|\gg X^{\be_n+\ep}$ for some $\be_n>0$, providing the first  lower bound that grows like a power of $X$.
Here is an open problem: for many groups $G$, no such polynomial $f$ has yet been exhibited.

Now we focus on the conjectured upper bound (\ref{weak_Malle}) for counting fields with bounded discriminant.
For any family $\Fscr = \Fscr_n(G;X)$ of fields,   the strong ``pointwise'' property $\Dbf_{\Fscr,n}(0)$ implies  Malle's ``average'' upper bound (\ref{weak_Malle}) for the group $G$; see   \cite{KluWan20x}. 
What is more surprising is that there is a converse to this. This relates to our  question: can the method of moments be used to deduce  the Discriminant Multiplicity Conjecture? Formally, it can. Given a family $\Fscr$ of fields, for each integer $D \geq 1$ let $m(D)$ denote the number of fields $K\in \Fscr$ with $D_K=D.$ If arbitrarily high $k$-th moment bounds are known for the function $m(D)$, the Discriminant Multiplicity Conjecture follows; see \cite{PTBW21}. But the first moment  of $m(D)$ is the subject of the Malle Conjecture  (\ref{weak_Malle}), so the method of moments certainly seems a difficult avenue to pursue. Yet interestingly,
 Ellenberg and Venkatesh have shown that in this context, $k$-th moments can be repackaged as    \emph{averages}.
 
Informally, the idea is to replace bounding the $k$-th moment of the function $m(D)$ for $G$-Galois fields in a family $\Fscr$ by counting fields in a family  $\Fscr^{(k)}$ of $G^k$-Galois fields. Ellenberg and Venkatesh order  the  fields in $\Fscr^{(k)}$   not by  discriminant $D_K$, but (roughly speaking) by the square-free kernel $D_K^\#$ of the   discriminant. They generalize the Malle Conjecture to posit that in this ordering,  $\ll X^{1+\ep}$ fields in $\Fscr^{(k)}$ have $D_K^\# \leq X$,  \emph{uniformly} for all integers $k \geq 1.$ Assuming this conjecture, suppose   there are $m(D)$ many $G$-Galois  fields $K_1,\ldots,K_{m(D)}$   with $D_{K_i} = D$.  Taking composita of $k$ of these generates at least $ \gg_k m(D)^k$ many $G^k$-Galois fields in the family $\Fscr^{(k)}$, with  $D_K^\# \leq D$. If we suppose  $m(D) \geq D^{\al}$ for some $\al>0$ and a sequence of $D \maps \infty,$ under the generalized Malle Conjecture it must be that $\al k \leq 1$ for all $k \geq 1$. Hence $\al$ must be arbitrarily small, as desired.

In full generality, Ellenberg and Venkatesh propose a generalized Malle Conjecture in terms of  an $f$-discriminant, for any rational class function $f$, and an appropriate generalization $a_G(f)$ of the exponent in (\ref{weak_Malle}). They verify that for a particular choice of $f$, this implies the Discriminant Multiplicity Conjecture. More recently, Kl\"uners and Wang have shown directly that Malle's Conjecture (\ref{weak_Malle}) for all groups $G$ implies the Discriminant Multiplicity Conjecture (also over any number field) \cite{KluWan20x}.

Let us sum up:  the upper bound (\ref{weak_Malle}) in Malle's Conjecture for all groups $G$ implies the Discriminant Multiplicity Conjecture. The Discriminant Multiplicity Conjecture implies the $\ell$-torsion Conjecture.  Also, the Discriminant Multiplicity Conjecture for $\Fscr_n(G;X)$  (that is, property $\Dbf_{\Fscr,n}(0)$) implies Malle's Conjecture for $\Fscr_n(G;X)$. 
Moreover there is one more converse: Alberts has shown that if the $\ell$-torsion Conjecture is true for all solvable extensions and all primes $\ell$ (even just in an average sense), then  Malle's upper bound (\ref{weak_Malle}) holds for all solvable groups \cite{Alb20}. Thus Malle's Conjecture, the Discriminant Multiplicity Conjecture, and the $\ell$-torsion Conjecture are truly equivalent, when restricted to solvable groups.
  These relationships provide  clear motivation for why so many methods described in this survey have involved counting number fields.

In conclusion, we have seen from three different perspectives that   the $\ell$-torsion Conjecture should be true.   But as Gauss wrote, ``\emph{Demonstrationes autem}
 rigorosae  \emph{harum observationum perdifficiles esse videntur}.''

\subsection*{Funding}
This work was partially supported by NSF CAREER DMS-1652173.

%***************************************
\bibliographystyle{alpha}
\bibliography{NoThBibliography}
%***************************************

\end{document}